\documentclass[11pt]{article}
\usepackage[T1]{fontenc}
\usepackage{lmodern}
\usepackage[margin=1in]{geometry}
\usepackage{amsmath,amssymb,amsthm}
\usepackage{booktabs}
\usepackage{graphicx}
\usepackage{xcolor}
\usepackage{algorithm}
\usepackage{algpseudocode}
\usepackage[hidelinks]{hyperref}
\usepackage{url}

\newtheorem{theorem}{Theorem}

\newcommand{\R}{\mathrm{R}}
\newcommand{\sha}[1]{\texttt{\scriptsize #1}}
\newcommand{\wR}{\omega_{\mathrm R}}
\newcommand{\wB}{\omega_{\mathrm B}}

\title{Lower bounds for Ramsey numbers:\\
\texorpdfstring{$\R(6,8)\ge 135$}{R(6,8) >= 135} and \texorpdfstring{$\R(8,10)\ge 345$}{R(8,10) >= 345}}
\author{Fritz Cremer}
\date{}

\begin{document}
\maketitle

\begin{abstract}
We prove the lower bounds $\R(6,8)\ge 135$ and $\R(8,10)\ge 345$, improving
the bounds 134 and 343 listed in the April 2026 revision of Radziszowski's
dynamic survey. We give explicit red/blue colorings of $K_{134}$ with no red
$K_6$ or blue $K_8$, and of $K_{344}$ with no red $K_8$ or blue $K_{10}$,
and verify them with two independently written exhaustive clique checkers.
Starting from published colorings, we find these witnesses by adding vertices
and repairing the resulting monochromatic cliques through local search.
The search uses exact conflict counts, preparation moves aimed at making
remaining cliques cheaper to break, and mutations followed by repair. A loss
based on the largest monochromatic clique containing each edge yielded a useful
intermediate state for the $\R(6,8)$ construction. Guided by the author, an AI
coding agent wrote and ran the search code; we describe the methods and
document the ancestry of the resulting colorings.
\end{abstract}

\section{Introduction}

The Ramsey number $\R(r,s)$ is the smallest $n$ such that every coloring of the
edges of $K_n$ with red and blue contains a red $K_r$ or a blue $K_s$. A
coloring of $K_n$ with neither, an \emph{$(r,s)$-coloring}, shows
$\R(r,s)\ge n+1$.

\begin{theorem}\label{thm:main}
$\R(6,8)\ge 135$ and $\R(8,10)\ge 345$.
\end{theorem}

Revision 18 of Radziszowski's survey \cite{radziszowski2026ds1}, dated 24 April
2026, lists $\R(6,8)\ge 134$, due to Exoo and Tatarevic
\cite{exootatarevic2015}, and $\R(8,10)\ge 343$, due to Kuznetsov
\cite{kuznetsov2016}. The corresponding upper bounds listed there are 423 and
4402. We prove the theorem with one $(6,8)$-coloring of $K_{134}$ and two
$(8,10)$-colorings of $K_{344}$. Four $(8,10)$-colorings of $K_{343}$, found
earlier, establish the intermediate bound $\R(8,10)\ge 344$; one of them
served as the starting point for both colorings of $K_{344}$.
Section~\ref{sec:colorings} presents and verifies all seven colorings.

Our searches start from published colorings, add a vertex, and repair the
resulting monochromatic cliques (Sections~\ref{sec:methods}
and~\ref{sec:ancestry}). This approach follows a long history of computational
constructions \cite{exootatarevic2015, kuznetsov2016, exoo2023r56}.
More recently, Nagda, Raghavan and Thakurta
improved nine bounds of the form $\R(3,s)$ and $\R(4,s)$ with AlphaEvolve
\cite{nagda2026alphaevolve}. Here an AI coding agent, guided by the author,
wrote and ran the search code (Section~\ref{sec:contributions}).

\section{The colorings}\label{sec:colorings}

Each coloring of $K_n$ is stored as a symmetric $n\times n$ matrix with zero
diagonal and vertices $0,\dots,n-1$, where an off-diagonal $1$ means red and $0$
means blue. The files are in the public repository \cite{cremer2026repo},
together with a graph6 copy of each red graph. Appendix~\ref{app:hashes} lists
the SHA-256 hashes of the matrix files, and Table~\ref{tab:colorings} gives their
red edge counts and maximum clique sizes. Here $\wR$ and $\wB$ denote the sizes
of the largest red and blue cliques.

\begin{table}[ht]
\centering
\small
\begin{tabular}{@{}lllrcc@{}}
\toprule
Bound & & File & Red edges & $\wR$ & $\wB$ \\
\midrule
$\R(6,8)\ge 135$ & & \texttt{R6\_8/n134-gpt-6.1-sol.txt} & 3\,913 & 5 & 7 \\
\addlinespace
$\R(8,10)\ge 345$ & no.\ 1 & \texttt{R8\_10/n344-gpt-6.1-sol-1.txt} & 25\,854 & 7 & 9 \\
 & no.\ 2 & \texttt{R8\_10/n344-gpt-6.1-sol-2.txt} & 25\,855 & 7 & 9 \\
\addlinespace
$\R(8,10)\ge 344$ & no.\ 1 & \texttt{R8\_10/n343-gpt-6.1-sol.txt} & 25\,770 & 7 & 9 \\
 & no.\ 2 & \texttt{R8\_10/n343-gpt-6-pro-1.txt} & 25\,809 & 7 & 9 \\
 & no.\ 3 & \texttt{R8\_10/n343-gpt-6-pro-2.txt} & 25\,795 & 7 & 9 \\
 & no.\ 4 & \texttt{R8\_10/n343-gpt-6-pro-3.txt} & 25\,795 & 7 & 9 \\
\bottomrule
\end{tabular}

\caption{The colorings. $\wR$ and $\wB$ were computed with Cliquer.}
\label{tab:colorings}
\end{table}

The colorings with the same $n$ are pairwise non-isomorphic: they differ in the
number of red edges or, for $K_{343}$ nos.\ 3 and 4, in their red-degree
sequences. All seven red graphs also have nonconstant degrees, so none is
circulant. Coloring no.~1 of $K_{343}$ was found
by the search of Section~\ref{sec:ancestry}; nos.~2 to 4 came from separate
ChatGPT sessions (Section~\ref{sec:contributions}).

\paragraph{Verification.} Each coloring was checked by two independently
written exhaustive searches for a red $K_r$ and a blue $K_s$, one based on
bitset intersection and one on an adjacency matrix, with no code in common. As
an external check, Cliquer~1.21 \cite{ostergard2002cliquer,cliquer} computed
$\wR$ and $\wB$. Table~\ref{tab:verification} records these checks; the
repository also contains a Python verifier based on NetworkX
\cite{hagberg2008networkx}.

\begin{table}[ht]
\centering
\small
\begin{tabular}{@{}lrrrrrr@{}}
\toprule
& \multicolumn{2}{c}{\texttt{clique\_oracle}} & \multicolumn{2}{c}{\texttt{verify\_independent}} & \multicolumn{2}{c}{Cliquer (s)} \\
\cmidrule(lr){2-3}\cmidrule(lr){4-5}\cmidrule(l){6-7}
Coloring & nodes & s & nodes & s & red & blue \\
\midrule
$(6,8)$, $n=134$ & 631\,956 & 0.030 & 91\,758 & 0.060 & 0.005 & 0.018 \\
$(8,10)$, $n=344$, no.\ 1 & 50\,255\,350 & 4.8 & 7\,030\,180 & 6.2 & 0.34 & 7.8 \\
$(8,10)$, $n=344$, no.\ 2 & 50\,227\,192 & 4.6 & 6\,690\,114 & 5.7 & 0.32 & 7.9 \\
$(8,10)$, $n=343$, no.\ 1 & 48\,333\,365 & 4.4 & 6\,982\,843 & 5.9 & 0.34 & 7.3 \\
$(8,10)$, $n=343$, no.\ 2 & 48\,189\,579 & 4.4 & 7\,218\,446 & 6.2 & 0.33 & 7.2 \\
$(8,10)$, $n=343$, no.\ 3 & 48\,215\,638 & 4.5 & 7\,180\,124 & 6.1 & 0.33 & 7.2 \\
$(8,10)$, $n=343$, no.\ 4 & 47\,773\,459 & 4.4 & 6\,886\,698 & 5.9 & 0.33 & 7.4 \\
\bottomrule
\end{tabular}

\caption{Verification of the public files on 6 October 2026, on an AMD Ryzen 9
5900X with g++ 13.3.0 (\texttt{-O3 -std=c++17}). ``Nodes'' are search-tree
nodes; times are single-run wall-clock seconds. All checks succeeded.}
\label{tab:verification}
\end{table}

\paragraph{Other pairs.} Applying monotonicity and standard lower-bound
inequalities from Section~2.3 of the survey, we found no further improved bounds
for $r,s\le 10$.

\section{Search methods}\label{sec:methods}

The searches allowed every edge to change color and imposed no symmetry or
fixed subgraph. We call a red $K_r$ or a blue $K_s$ a \emph{conflict}.
Flipping an edge $e$ changes the numbers of red and blue
conflicts by $\Delta(e)=(\Delta_{\mathrm R}(e),\Delta_{\mathrm B}(e))$; for
example, flipping a red edge $uv$ removes every red $K_r$ through it and creates
one blue $K_s$ for each blue $K_{s-2}$ among the common blue neighbours of $u$
and $v$. The descriptions below follow the source code of the programs that
ran.

\subsection{Exact conflict repair}

The basic engine maintains the complete sets of red and blue conflicts and
updates them after each flip, giving exact deltas. Initially
$w_{\mathrm R}=s^2/r^2$ and $w_{\mathrm B}=1$ to account for the different clique
sizes. Candidate edges are sampled from current conflicts, although a flip may
create new conflicts elsewhere. A Metropolis rule permits increases in the
weighted conflict count, with a temperature that rises during stalls.
Algorithm~\ref{alg:repair} specifies one iteration with sampling budget $m$ and
stall threshold $k$; an ordinary iteration is skipped if all sampled edges are
tabu and fail the aspiration rule.

\begin{algorithm}[ht]
\caption{One iteration of exact conflict repair.}
\label{alg:repair}
\begin{algorithmic}[1]
\State let $t$ count iterations since the last improvement or kick
\If{$t>k$}
  \State make up to 8 conflict-directed random flips, clear tabu memory, reset $t$, and set $w_{\mathrm R}=(s^2/r^2)e^U$ with $U$ uniform on $[-0.8,0.8]$; end this iteration
\EndIf
\State choose a color $c$ with probability proportional to $w_c\,|\mathcal K_c|$, where $\mathcal K_c$ is the set of conflicts of color $c$, and a conflict $K\in\mathcal K_c$ uniformly at random
\State sample up to $m$ edges of $K$ and compute $\delta(e)=w_{\mathrm R}\Delta_{\mathrm R}(e)+\Delta_{\mathrm B}(e)$ for each, skipping tabu edges unless they beat the best count so far
\State let $e^\ast$ minimize $\delta$ and set $T=0.25+\min(2.5,t/100)$
\If{$\delta(e^\ast)\le 0$ \textbf{or} a uniform random number is below $\exp(-\delta(e^\ast)/T)$}
  \State flip $e^\ast$ and make it tabu for 5 to 19 iterations
\EndIf
\end{algorithmic}
\end{algorithm}

\subsection{Preparation}

Near a solution, breaking one conflict often creates several others. We call
the minimum net change in conflict count,
$\min_e (\Delta_{\mathrm R}(e)+\Delta_{\mathrm B}(e))$ over edges $e$ of current
conflicts, the \emph{breaking cost}. The following methods aim to make these
conflicts cheaper to break.

\emph{Global preparation} runs periodically inside Algorithm~\ref{alg:repair}.
It samples up to 32 edges of current conflicts and chooses $uv$ minimizing the
weighted delta $\delta$. If this minimum is positive, the kernel lists the
cliques that flipping $uv$ would create and tries flips of other edges that meet
them. Each preparation flip preserves both conflict sets and is kept if it
lowers the minimum $\delta$ over the sampled edges; ties are kept with probability
$1/4$. After at most 8 such steps, the cheapest edge in the updated sample is
proposed under the acceptance rule of Algorithm~\ref{alg:repair}.

A \emph{risk chain} has the same goal but uses up to three flips that may
temporarily create up to two extra conflicts. It is committed only if the final
set of conflicts is exactly the original one and the cost of breaking the
target edge has strictly decreased.

\emph{Clique weighting with ejection chains}, used early in the $\R(6,8)$ search,
scores flips by individual conflict weights. These weights are periodically
increased for surviving conflicts, with a cap and occasional decay, as in
breakout methods. An ejection chain first flips an edge of a conflict, then
repeatedly flips an edge of a conflict created by the previous flip. The chain
is kept only if the number or total weight of conflicts drops.

\subsection{Mutation followed by repair}

When repair stalls, a mutation changes a small region before unrestricted
repair resumes. This adapts the section mutation and crossover followed by
relaxation in Jeroen Cottaar's winning solution to the Santa 2025 Kaggle
competition \cite{cottaar2025santa}; no code from it was used.
For section mutations, a set $S$ of up to 8 vertices is built around a vertex
chosen with probability proportional to its number of conflicts. A
\emph{section crossover} copies donor colors from the same lineage on up to a
given number of pairs meeting $S$, while a \emph{patch} flips random pairs
inside $S$. Two further moves act on selected vertices: a \emph{rectangle}
flips red $ab,cd$ and blue $ac,bd$ together, preserving all red degrees;
a \emph{block move} exchanges the adjacencies of two vertices to a selected
subset of other vertices.

\subsection{The max-clique edge loss}\label{sec:loss}

Conflict counts do not directly measure smaller cliques that may grow into new
conflicts. We therefore propose a loss in which each edge contributes the size
of the largest monochromatic clique containing it,
\[
  L=\sum_{uv}\omega(uv),\qquad
  \omega(uv)=2+\omega_c\bigl(N_c(u)\cap N_c(v)\bigr),
\]
where the sum is over unordered edges of $K_n$, $c$ is the color of $uv$,
$N_c(v)$ is the set of vertices joined to $v$ in color $c$, and $\omega_c(S)$
is the size of the largest clique of color $c$ inside $S$. The loss does not depend on
$(r,s)$. It is computed exactly with a branch-and-bound
clique routine, rescoring after each flip only the edges whose value can change.
In the loss-guided run used in the $(6,8)$ construction, conflict edges were
sampled with probability $0.8$; otherwise, background moves preserving both
conflict counts were sampled. Candidates were ranked by their change $\Delta L$
and accepted with Metropolis probability $\min(1,e^{-\Delta L/T})$.
Here $T=0.5+\min(8,t/100)$, with $t$ as in Algorithm~\ref{alg:repair}; conflict
counts still governed tabu aspiration, stall detection and checkpoint selection.

In short controlled comparisons from a seven-conflict $(6,8)$ state, neither
count-only nor loss-guided search reduced the conflict count. Using the loss as
the primary objective made search eight to nine times slower; using it only to
break count ties had smaller overhead. In that starting state, almost every
edge lay in a near-forbidden clique, so the loss largely tracked edge density.
One of these runs supplied the intermediate state used in the $(6,8)$
construction described in Section~\ref{sec:r68}.

\section{How the colorings were found}\label{sec:ancestry}

For each saved search checkpoint, we retained its command, random seed and
source code. Figures~\ref{fig:r810} and~\ref{fig:r68} trace the resulting
colorings back to their published starting points. Each link from a coloring
to its successor is recorded by file hash, except for twelve handovers between
early search segments, which are identified from the search supervisor's log.
Saved move traces allow some trajectories to be replayed,
including the 43-flip loss-guided run. Exact reruns of the whole campaign are
not guaranteed: runs were limited by wall-clock time and full internal search
states were not checkpointed. The campaigns varied in seeds, duration and
parents, so their histories document successful routes without providing a
controlled comparison of search methods.

\subsection{\texorpdfstring{$\R(8,10)$}{R(8,10)}}\label{sec:r810}

\begin{figure}[!htbp]
\centering
\includegraphics[width=\textwidth]{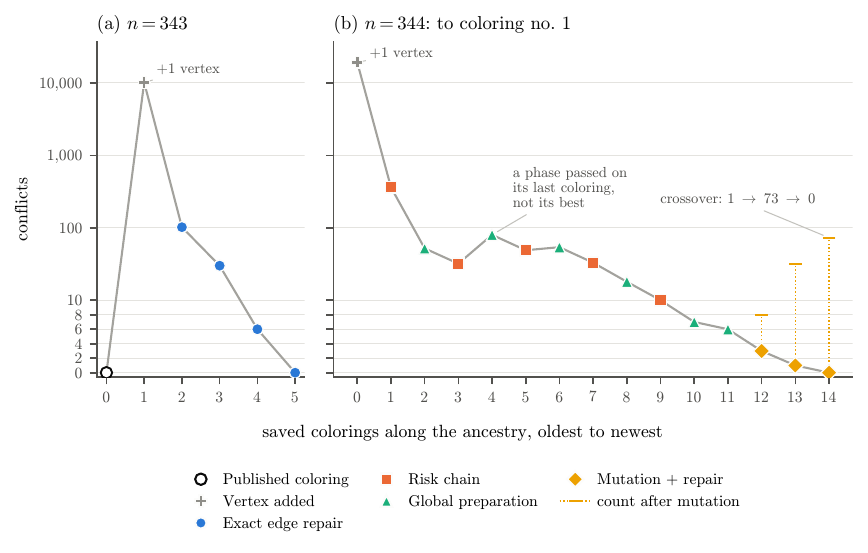}
\caption{Ancestry of the $(8,10)$-colorings: (a) from Kuznetsov's coloring of
$K_{342}$ to coloring no.~1 of $K_{343}$; (b) from that coloring plus one vertex
to the seven-conflict parent of both colorings of $K_{344}$ (step 10), and on
to coloring no.~1; the branch to no.~2 is similar and not shown. Each marker is a saved coloring, colored and shaped by the method
that produced it; the vertical axis is the exact number of conflicts, linear up
to 10 and logarithmic above. A dotted orange line rises from a coloring to
the count right after the mutation that preceded it; repair then brought the
count down to the marker. Phases that returned their input unchanged are
omitted.}
\label{fig:r810}
\end{figure}

Kuznetsov's $(8,10)$-coloring of $K_{342}$ \cite{kuznetsov2016} is circulant:
$\{i,j\}$ is red exactly when $\min(|i-j|,342-|i-j|)$ lies in a set of 76
distances. Adding a copy of vertex 0, joined to it by a blue edge, gives 10\,146
conflicts. Four consecutive runs of exact repair, about half an hour of search
in total, reduced this to 102, 30, 6 and 0 (Figure~\ref{fig:r810}(a)). The
result, coloring no.~1 of $K_{343}$, differs from the initial copy on 874 edges.

We extended this coloring by copying vertex 269 and joining the two copies by
a blue edge. The resulting coloring of $K_{344}$ has 19\,271 blue conflicts.
Ten 30-minute segments alternating risk chains and global preparation reduced
the count to 7 (Figure~\ref{fig:r810}(b)). From this shared parent, two
global-preparation runs saved states with 6 and 2 conflicts, respectively.
Further mutation and repair led to the two witnesses of $K_{344}$. For no.~1
(Figure~\ref{fig:r810}(b)), three section crossovers with repair,
each using the seven-conflict coloring as donor, gave $6\to3\to1\to0$; the last
crossover copied 26 donor colors and briefly raised the count from 1 to 73
before repair reached zero. For no.~2, two patches with repair gave
$2\to1\to0$; the last patch flipped four edges and briefly raised the count
from 1 to 140.

\subsection{\texorpdfstring{$\R(6,8)$}{R(6,8)}}\label{sec:r68}

\begin{figure}[!htbp]
\centering
\includegraphics[width=\textwidth]{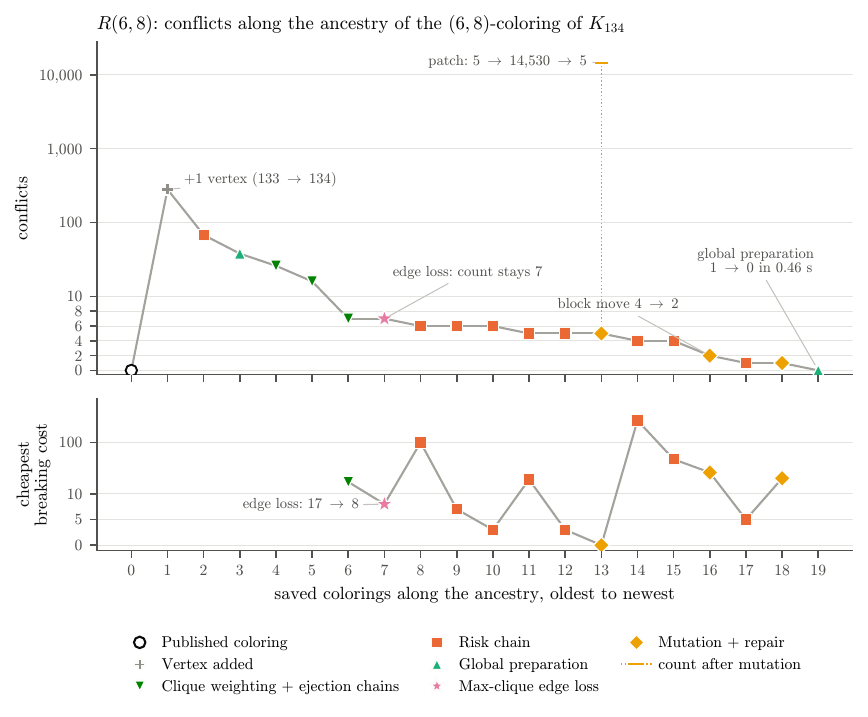}
\caption{Ancestry of the $(6,8)$-coloring, with the conventions of
Figure~\ref{fig:r810}. Top: number of conflicts. Bottom: breaking cost, where
recorded. The max-clique edge loss step (star) left the count at 7 but lowered
the breaking cost from 17 to 8.}
\label{fig:r68}
\end{figure}

The search started from the $(6,8)$-coloring of $K_{133}$ by Exoo and Tatarevic
\cite{exootatarevic2015}. Adding a copy of vertex 117, joined to it by a red
edge, gives 281 red conflicts. A 30-minute risk-chain segment reduced the total
conflict count to 68, a global-preparation segment to 38, and three segments of
clique weighting with ejection chains to 26, 16 and 7 (Figure~\ref{fig:r68}).

From 12 endpoints that still had 7 conflicts, we selected the one with the
lowest breaking cost. It came from a run guided by the max-clique edge loss
that made 43 flips, keeping the same 7 conflicts while lowering the breaking
cost from 17 to 8. Two 300-second risk-chain runs from this state both reached
6 conflicts. Continued preparation then lowered the breaking cost of the
six-conflict state from 100 to 7.

The prepared six-conflict coloring seeded a campaign of six parallel lanes,
each switching move family every 170 seconds and sharing its best colorings, ranked by conflict
count and then by breaking cost. Along the ancestry, risk-chain phases reduced
the count from 6 to 5 to 4. A block-exchange repair phase reached 2, followed by
a risk-chain phase that reached 1. A final global-preparation phase eliminated
the last conflict in 0.46 seconds.

\section{Human and AI contributions}\label{sec:contributions}

The search code was generated and run by Codex, OpenAI's coding agent, using the
model GPT-6.1 Sol with reasoning effort ``Extra High''. The author guided the
search, setting the objectives and proposing ideas and strategies,
including the transfer of section crossover from the Kaggle solution and the
max-clique edge loss. The agent ran on the author's AMD Ryzen 9 5900X
workstation and later on up to 76 CPU cores at a time.

Coloring no.~1 of $K_{343}$ was found less than two hours after the agent was
first asked to improve a bound, and nos.~2 to 4 were produced by GPT-6 Pro in
separate ChatGPT sessions from the single prompt ``Improve the Ramsey lower bound
for R(8, 10)''.

AI tools were used as writing and editing aids during preparation of the
manuscript.

\section*{Data availability}

The colorings, as matrices and in graph6 format, and a Python verifier are in
the repository \url{https://github.com/fritzcremer/ramsey-lower-bounds}, release
\texttt{v1.1.0} \cite{cremer2026repo}.

\appendix

\section{File hashes}\label{app:hashes}

The SHA-256 hashes of the matrix files in release \texttt{v1.1.0} of the
repository \cite{cremer2026repo} are as follows.

\begin{table}[ht]
\centering
\small
\begin{tabular}{@{}ll@{}}
File & SHA-256 \\
\midrule
\texttt{n134-gpt-6.1-sol.txt} & \sha{8ef378ad4d66478c33939e438d3a64572ec2cf7d1647ae29d3568817fedfc982} \\
\texttt{n344-gpt-6.1-sol-1.txt} & \sha{17881b78e3d8b55372995f4e53ac068ca22560cf30b291c1e649eacaef21c21d} \\
\texttt{n344-gpt-6.1-sol-2.txt} & \sha{fd3c5a55434179a596dacd9f94a1208211e3ce0da7eee9502c640a3b1860ad89} \\
\texttt{n343-gpt-6.1-sol.txt} & \sha{f53fa7d45f427f16a159660194cf2de5b875fa0153660eb0fe348984e450095b} \\
\texttt{n343-gpt-6-pro-1.txt} & \sha{b913b14121c880a55f957947451cb42bcf1cf02ee3bc4fb40a6eb44d7661e0e1} \\
\texttt{n343-gpt-6-pro-2.txt} & \sha{4c0c0b7b4df4cd4040e20bca7c724a790e9ed35ae05aafdc30f2192806e840ec} \\
\texttt{n343-gpt-6-pro-3.txt} & \sha{a762539bb9d8543fb1313139c9b87e3e5ee04530aaa2faffa84aa5b9ed6c4263} \\
\end{tabular}
\end{table}

\bibliographystyle{plain}
\bibliography{refs}

\end{document}